\documentstyle[12pt]{report} 

\newtheorem{lemma}{Lemma}
\def\blm{\begin{lemma}}
\def\elm{\end{lemma}}
\newtheorem{theorem}{Theorem}

\def\beq {\begin{equation}}
\def\eeq {\end{equation}}
\def\btm{\begin{theorem}}
\def\etm{\end{theorem}}
\def\ben{\begin{enumerate}}
\def\een{\end{enumerate}}

\def \ct#1 {\mbox{{$Cob_3($#1$)^{conn}\,\,$}}}

\input{epsf.sty}

\begin{document}

\ \hfill{{{\small To appear in } {\em Topology Appl.}}}

\ 

\begin{center}

\section*{Equivalence of a Bridged Link Calculus\\ and Kirby's Calculus of
Links\\   
on Non-Simply Connected 3-Manifolds}

\vspace*{.5cm}
\bigskip

{\large Thomas Kerler}\\
\medskip

 November 1994 (Revised Version April 1997)

\end{center}
 \vspace*{.5cm}

{\small \noindent{\bf Abstract} : We recall an extension of Kirby's
Calculus on non-simply connected 3-manifolds given in [FR], and the
 surgery calculus of bridged links from [Ke], which involves only local 
moves. We give a short combinatorial proof that the two calculi
are equivalent,  and thus describe the same classes of 3-manifolds.
This makes the proofs for the validity of surgery calculi in [FR] and [Ke]
interchangeable.  }

\paragraph*{1.) \underline{The Extended Kirby Calculus} :}
 
A popular  way of presenting a closed, compact, connected, and oriented 
3-manifold, $\, N^{(3)}\,$,
is by use of framed links in $\,S^3\,$. To
any such link, $\,L\,$, we can associate a three-fold, $\,S^3_L\,$,
by doing {\em surgery} along its components. It is long been known
[LW] that  any  closed $\, N^{(3)}\,$ is homeomorphic to 
 $\,S^3_L\,$ for some link, $\,L\,$. 

 In [Ki] Kirby shows  that $\,S^3_{L_1}\,$ is homeomorphic to 
$\,S^3_{L_2}\,$, if and only if $\,L_1\,$ can be changed into  $\,L_2\,$
by a sequence of certain {\em moves}, i.e., if   $\,L_1\,$  and  $\,L_2\,$ 
are {\em $\delta$-equivalent } in the language of [Ki]. The moves are either
ambient isotopies, the {\em signature  move} ${\cal O}_1$, or the {\em two-handle
slide} $\,{\cal O}_2\,$. Thickening the links into ribbons, in order to
indicate their framings, the $\,{\cal O}_2\,$-move is depicted in Figure 1, 
where $\,R_2\,$ is slid over $\,R_1\,$.

\begin{center}\ 
\epsfbox{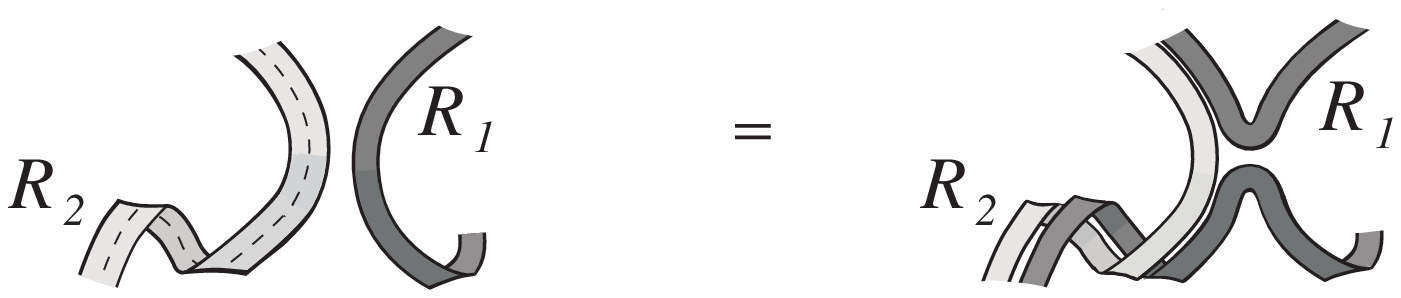}
\newline{{\sc Figure 1:} $\,{\cal O}_2\,$-{\em move}}
\end{center}

Hence Kirby's calculus establishes a one-to-one correspondence between the
set of homeomorphism classes of closed 3-manifolds and the combinatorial
set of  $\delta$-equivalence classes of links in $\,S^3\,$.
\smallskip

In a more general situation we can consider framed links $\,L\,$ embedded into a connected, compact, and oriented 3-manifold, $\,M=M^{(3)}\,$, which may have a boundary or be non-simply connected. For  two 
framed links, $L_1$ 
and $\,L_2\,$, in $\,M\,$, we can again consider the correspondingly surgered manifolds
$\,M_{L_1}\,$ and $\,M_{L_2}\,$. It was suggested by Kirby
and proven by Fenn and Rourke [FR] that $\,M_{L_1}\,\cong\,M_{L_2}\,$, if and only 
if $\,L_1\,$ is {\em $\alpha$-equivalent}
 to $\,L_2\,$, where the  $\alpha$-equivalences
are generated by isotopies, the  ${\cal O}_1\,$-,  ${\cal O}_2\,$-, and  
${\cal O}_3\,$-move. The new ${\cal O}_3$-move (also called $\eta$-move in [Ke])
is given by including or
deleting a two-component configuration from a link, as shown in Figure 2.
It consist of an arbitrary component, $\,R\,$, together with an annulus, 
$\,A\,$,
that encircles $\,R\,$.

\begin{center}\ 
\epsfbox{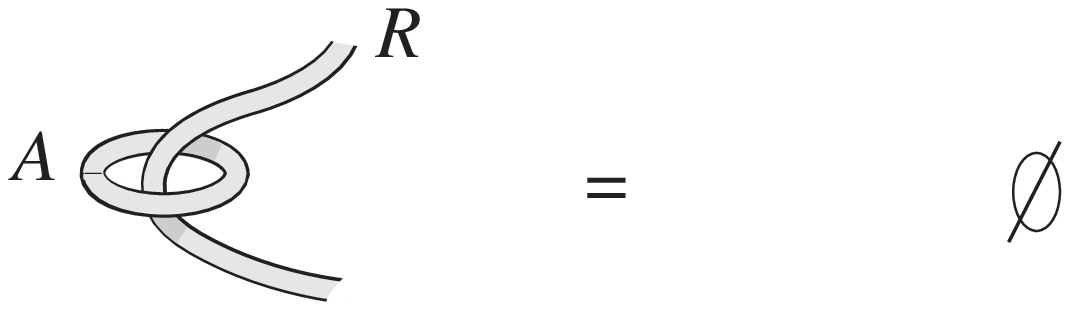}
\newline{{\sc Figure 2:} $\,{\cal O}_3\,$-{\em move}}
\end{center}

For the purpose of this letter let us consider the slightly smaller
equivalent classes, where we   omit the ${\cal O}_1$-move as a generator of
equivalences. The space of classes of framed links in a given
compact manifold, $\,M\,$, modulo  ${\cal O}_2$- and ${\cal O}_3$-moves
shall be  denoted $\,{\cal L}_M\,$.

It can be taken from  arguments in [Ki], [FR], and [Ke] that, in fact, it
suffices to choose a much smaller set of
  ${\cal O}_3$-moves by selecting  two
long ribbons $\,R\,$ for each class in $\,[S^1, M]\,$, with different framings
${\rm mod}\, 2\,$. All other ${\cal O}_3$-moves can be reduced to these special ones 
by combining them with isotopies and ${\cal O}_2$-moves. 
The ${\cal O}_3$-move for the trivial element in $\,[S^1, M]\,$ consists then of adding or deleting isolated Hopf links, $^0\bigcirc\mkern-13mu\bigcirc^{0,1}\,$, with framings 
as indicated.  Applying ${\cal O}_2$-moves they can both be substituted for the 
pair of  ${\cal O}_1$-moves, in which  $\bigcirc^{1}\sqcup\bigcirc^{-1}$ is 
added or deleted.

In order to see the topological interpretation of $\,{\cal L}_M\,$, recall
that a framed link $\,L\,\subset\,M\,$ defines also a four dimensional 
cobordism $\,W_L:M\to M_L\,$, with $\partial W_L\cong M \cup_{\partial M}M_L\,$.
It is now easily extracted from [Ki] or [FR] that $\,L_1\,$ and  $\,L_2\,$
are equivalent via  ${\cal O}_2\,$-, and  
${\cal O}_3\,$-moves, if and only if $\,M_{L_1}\cong M_{L_2}\,$ {\em and}
$\,{\rm sign}(W_{L_1})={\rm sign}(W_{L_2})\,$ for the signatures. The 
${\cal O}_1$-move connects a ${\bf CP}^2$ to $\,W_L\,$ and thus increases 
${\rm sign}(W_L)\,$ by one. (Note, however, that
 connecting ${\bf CP}^2\#\overline{{\bf CP}}^2$ 
or $S^2\times S^2\,$, as in 
the $^0\bigcirc\mkern-13mu\bigcirc^{0,1}$-moves, does not change the signature).

 Moreover, if  $\,N\,$ is any other compact, oriented, and connected
3-manifold with $\,\partial M\cong \partial N\,$, then there is a framed link 
$\,L\,$, such that $\,M_L\cong N\,$, extending the given homeomorphism on the
boundaries. Hence the classes in $\,{\cal L}_M\,$ are 
in one-to-one correspondence
with pairs $\,(N,n)\,$, where $\,N\,$ is a homeomorphism class of 3-manifolds
with boundary $\,\cong\partial M\,$, and $n\in{\bf Z}\,$ is given by $\,n={\rm sign}(W_L)\,$.
Let us denote the set of all such pairs by $\,{\cal M}_{\partial M}\,$.  
\smallskip

The proof in [FR], although referring to results from [Ki], is quite involved since
$\alpha$-equivalence of two links is first expressed as  $\delta$-equivalence 
together with an isomorphism between the corresponding $\pi_1(W_L)\,$ and
the vanishing of an additional obstruction class. A shorter proof that builds
directly on the one in  [Ki] was given very recently by Roberts [R]. Several
central arguments used there to generalize the steps in the proof of [Ki] had 
already  been given  in [Ke] and an earlier version of this letter.

\paragraph*{2.) \underline{Calculus of Bridged Links} :} 
In [Ke] we derive an alternative calculus of links, starting in a
similar way as in [Ki] but following a different strategy to single out
moves that generate the necessary equivalences. In order to include
the possibility of index-1-surgery we consider 
{\em bridged links} $\,L\subset M\,$ in a connected, compact, and oriented
3-manifold,  $M\,$. Thus $\,L\,$ not only consists of ribbons but also pairs
of balls, and a ribbon can enter one ball of  a pair and reemerge at
a corresponding spot at the other ball. 

In contrast to the non-local moves ${\cal O}_2$ and ${\cal O}_3$ of the 
extended Kirby calculus, the inclusion of 1-handle data allows us to 
use instead two {\em local} moves. The motivation in [Ke] is that 
locality is essential, e.g.,  for efficient presentations of cobordism categories, 
Hopf-algebraic interpretations of elementary cobordisms, and  thus 
constructions of TQFT's
by means of straight forward structural comparisons
 rather than involved calculations, see [KL].

  The easier of the two moves is the {\em isolated cancellation move}, which we
shall call here ${\cal A}_1\,$. As depicted in Figure 3 the 
${\cal A}_1\,$-move adds or deletes an isolated configuration from
the rest of the link, which consists of one ribbon, and one pair of
balls. The ribbon passes through the surgery spheres exactly once.

\begin{center}\ 
\epsfbox{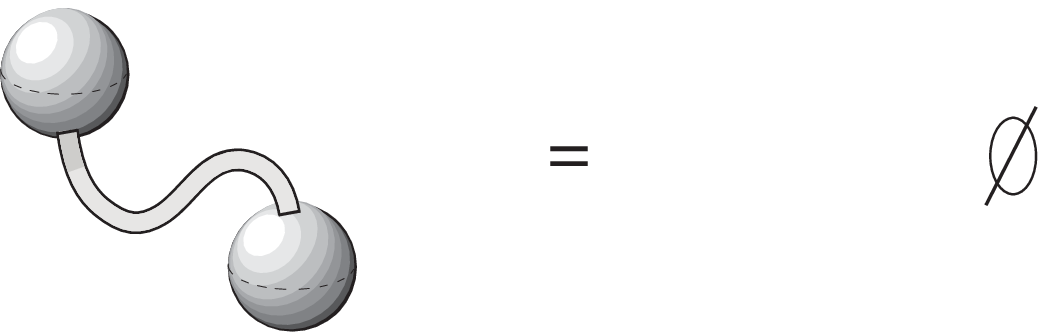}
\newline{{\sc Figure 3:} $\,{\cal A}_1\,$-{\em move}}
\end{center}

The ${\cal A}_2\,$-move, more commonly called {\em modification} or {\em handle trading}, is
described as follows. Any 0-framed component, $\,A\,$, bounding a disc,
can be replaced by a pair of surgery balls, $B_1$ and $B_2$, as indicated
in Figure 4. Other strands of the links that run through the disc bounded
by $\,A\,$ are redirected through the surgery spheres.

\begin{center}\ 
\epsfysize=2.8in
\epsfbox{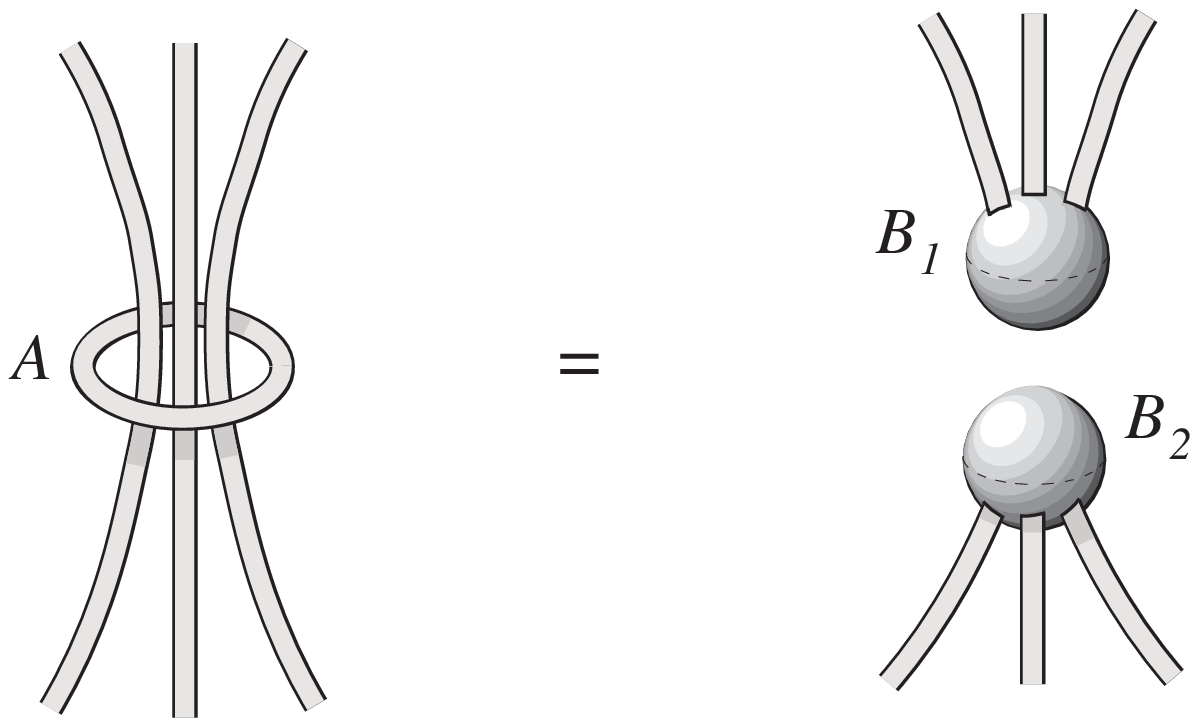}
\newline{{\sc Figure 4:} $\,{\cal A}_2\,$-{\em move}}
\end{center}

Let us denote by $\,{\cal BL}_M\,$ the space of bridged links modulo equivalences 
generated by isotopies, 
$\,{\cal A}_1\,$-, and $\,{\cal A}_2\,$-moves.
\medskip

 In [Ke] it was shown that $\,{\cal BL}_M\,$ for a connected, compact, and 
oriented 3-manifold is, like $\,{\cal L}_M\,$, in 
one-to-one correspondence with the set $\,{\cal M}_{\partial M}\,$ of pairs $\,(N,n)\,$, 
where $\,N\,$ is a homeomorphism class of 3-manifolds with 
$\partial N\cong\partial M\,$ and 
$\,n\in{\bf Z}\,$. 

The correspondence is again given by the assignment
$\,L\,\mapsto\,(M_L,{\rm sign}(W_L))\,$, only now we obtain $\,M_L\,$
by first performing index-1-surgeries along the pairs of balls
(yielding $\breve M\cong \,M\# S^1\times S^2\#\ldots\# S^1\times S^2\,$) before we do
the usual index-2-surgeries along the remaining framed link in $\,\breve M\,$. 
Similarly, $\,W_L\,$ is obtained, by also adding 1-handles instead of
just 2-handles to $\,M\times [0,1]\,$.

Since both link spaces have the same topological interpretation in
terms  of $\,{\cal M}_{\partial M}\,$ it is obvious that there has to be a one-to-one
correspondence also between $\,{\cal L}_M\,$ and $\,{\cal BL}_M\,$ for
a given connected $\,M\,$ as above. 

The main purpose of this letter is to present a direct and concise of
proof of the equivalence of  $\,{\cal L}_M\,$ and $\,{\cal BL}_M\,$ as
{\em formal} link calculi in $\,M\,$, and thereby show, how the two calculi
can be translated into each other  purely combinatorially. More precisely,
we prove the following:

\begin{theorem}\ 

The identification of a framed link in $\,M\,$ as a special bridged link
 in $\,M\,$ factors into a map on the equivalence classes of links
$$
I\,:\;{\cal L}_M\,\longrightarrow\,{\cal BL}_M\;.
$$
Moreover, the map $\,I\,$ is a bijection. 
\end{theorem}

 In the proof of Theorem 1 we shall only use the definition of
$\,{\cal L}_M\,$ and $\,{\cal BL}_M\,$ as sets of classes of
links as given in Paragraphs 1 and 2. In particular,   
none of the arguments in the proof will involve the topological
interpretation of 
the links as  surgery data. 

In the next theorem let us relate this combinatorial result to the 
results  
on presentations of 3-manifolds outlined previously. The first part 
 is obvious from the fact that the moves 
${\cal O}_2\,$,  ${\cal O}_3\,$, ${\cal A}_1\,$, and ${\cal A}_2\,$
can be interpreted as manipulations of $\,W_L\,$, which do not change 
the homeomorphisms type of the boundary piece  $\,M_L\,$.

\begin{theorem}\ 
Let $\,M\,$ be a compact, connected, and oriented 3-manifold.
\begin{enumerate}
\item The interpretations of (bridged) links as surgery prescriptions
factor into  maps:
$$
S_{\rm KFR}\,:\;{\cal L}_M\,\longrightarrow\,{\cal M}_{\partial M}\qquad\;
{\rm and }\qquad\;
S_{\rm BL}\,:\;{\cal BL}_M\,\longrightarrow\,{\cal M}_{\partial M}
$$
with\vspace*{-.9cm}

$$\; S_{\rm KFR}\;=\;S_{\rm BL}\circ I\quad.$$

\item Both  $\,S_{\rm KFR}\,$ and $\,S_{\rm BL}\,$ are bijections.
\end{enumerate}
\end{theorem}

The second part reveals an interesting application of Theorem 1 in that 
bijectivity of $\,I\,$ allows us to exchange the  proofs for the
two surgery calculi. I.e., we can use the proofs in [FR] or [Ki]
in order to establish the validity
of the bridged link calculus, or , alternatively, infer the extended
Kirby calculus starting from the results in [Ke].

\paragraph*{4.) \underline{Proof of Theorem 1} :\ }\

{\em 4.1) Existence of $\,I\,$: } In order to see that $\,I\,$ is well defined,
 we need to
show that if two links lie in the same class in  $\,{\cal L}_M\,$, then they are
also equivalent as links in $\,{\cal BL}_M\,$. For example, if two links differ
by an $\,{\cal O}_3\,$-move, we see from Figure 5 that they can also be 
related by  a combination of  $\,{\cal A}_1\,$- and $\,{\cal A}_2\,$-moves,
and thus belong to the same class in $\,{\cal BL}_M\,$.

\bigskip

\begin{center}\ 
\epsfxsize=5.8in
\epsfbox{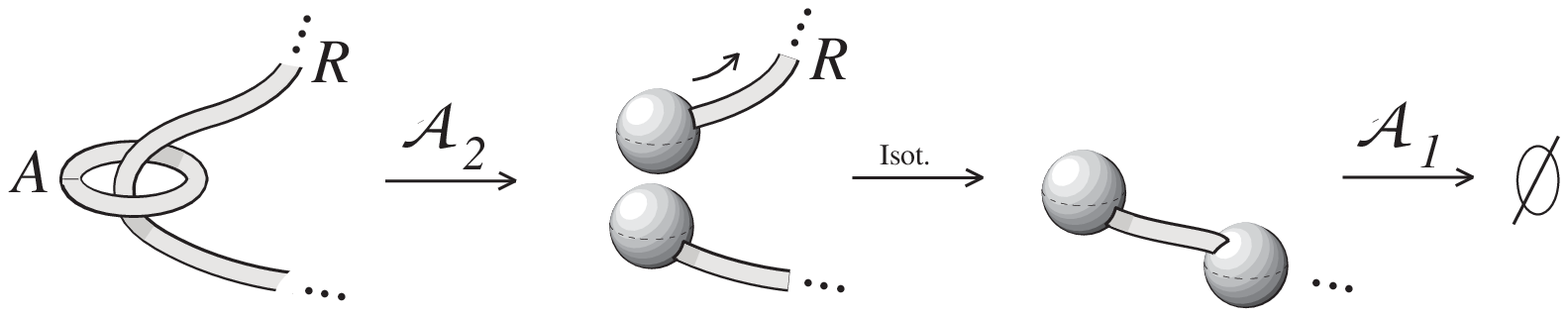}
\newline{{\sc Figure 5:}  $\,{\cal O}_3\,$ {\em from} $\,{\cal A}_1\,$ {\em and} 
$\,{\cal A}_2\,$ }
\end{center}

A more general version, $\,{\cal A}^{gen}_1\,$, of the cancellation move depicted in
Figure 3 is defined by letting additional strands pass through the surgery balls.
It is shown in Figure 6 that, with the $\,{\cal O}_3\,$-move now found to be an 
equivalence in $\,{\cal BL}_M\,$, also $\,{\cal A}^{gen}_1\,$ is an equivalence in
$\,{\cal BL}_M\,$.

\bigskip

\begin{center}\ 
\epsfxsize=5.8in
\epsfbox{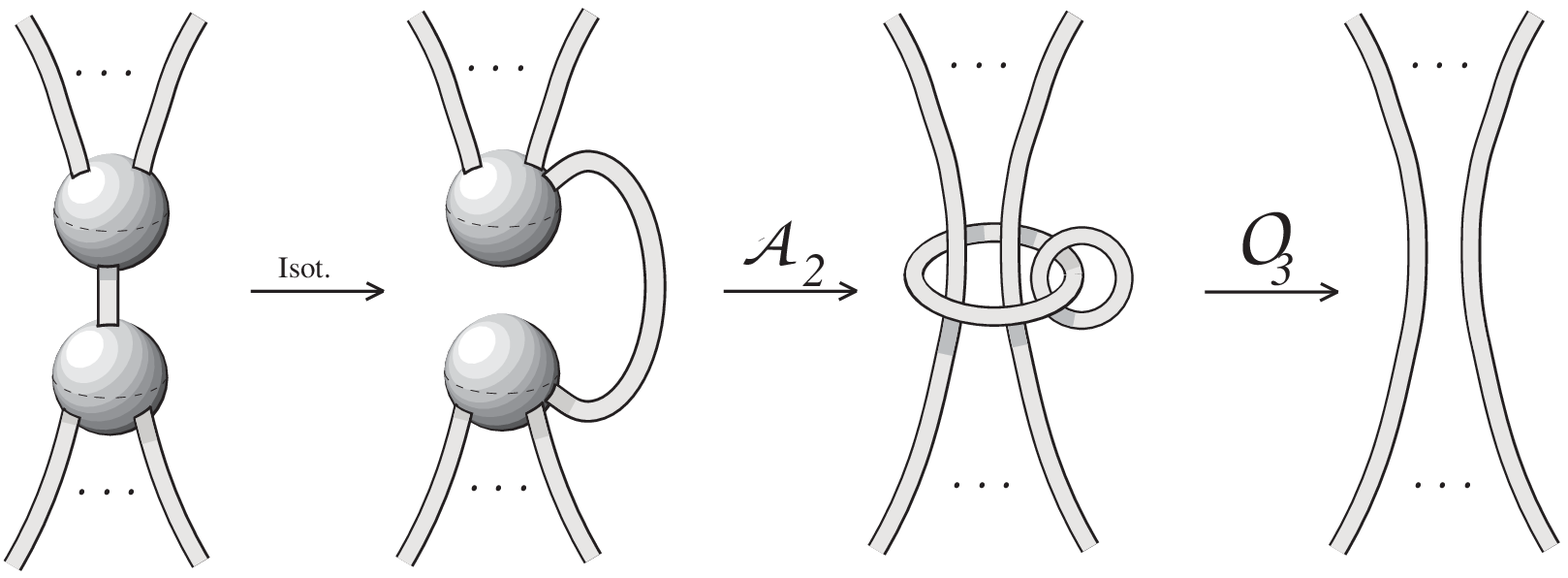}
\newline{{\sc Figure 6:} $\,{\cal A}_1^{gen}\,$ {\em from} $\,{\cal O}_3\,$ {\em and} 
$\,{\cal A}_2\,$ }
\end{center}

Moreover, as depicted in Figure 7,
the $\,{\cal O}_2\,$-move can be expressed as a combination of 
a general cancellation and uncancellation $\,{\cal A}^{gen}_1\,$, and 
an isotopy, where the surgery ball $\,B\,$ is pushed all the way along
the  second ribbon until it meets $\,B'\,$ again at the other end.

\bigskip

{\em 4.2) Surjectivity of $\,I\,$: } It is  easy to see that every 
class in $\,{\cal BL}_M\,$ has an ordinary, framed link as a representative.
Given an arbitrary bridged link, such a representative in the same class
can be found by applying an $\,{\cal A}_2\,$-move to each pair of balls after
they have been moved together in the connected manifold $\,M\,$. 
The choice of isotopy, by which the balls are moved
together, is conveniently expressed by a {\em recombination band} between
the two balls in the original bridged link. See for example Figure 8, where 
we indicated the isotopy and $\,{\cal A}_2\,$-move by the recombination band  
$\,r^B\,$ between $\,B\,$ and  $\,B'\,$.

\begin{center}\ 
\epsfxsize=6.1in
\epsfbox{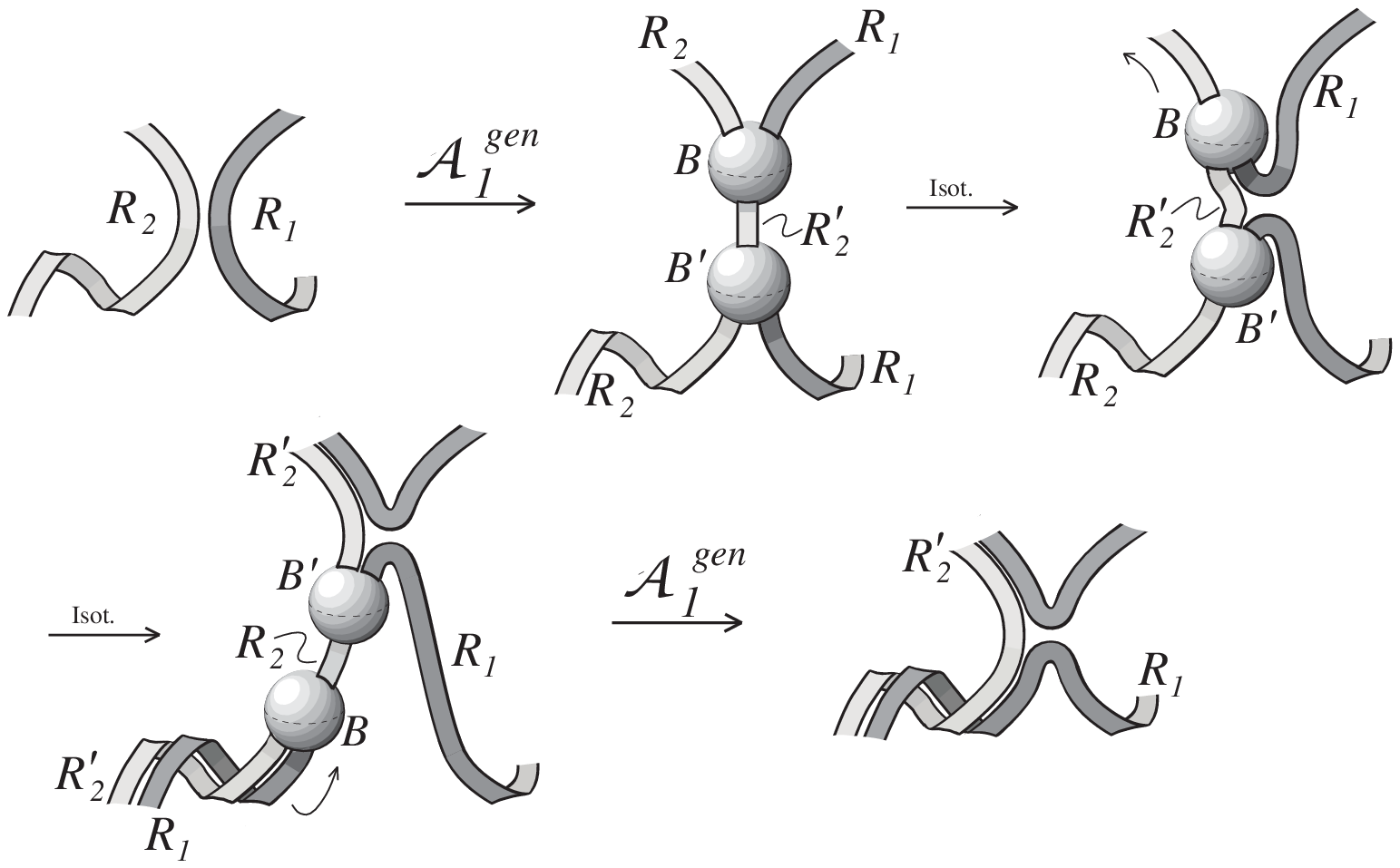}
\newline{{\sc Figure 7:} $\,{\cal O}_2\,$ {\em from} $\,{\cal A}_1^{gen}\,$}
\end{center}

{\em 4.3) Injectivity of $\,I\,$: } We verify injectivity by constructing  
a left inverse, $\,P\,$,  for $\,I\,$.  The image of $\,P\,$ of a class in
$\,{\cal BL}_M\,$ shall be defined by taking a representing bridged link,
recombining all pairs of surgery balls via $\,{\cal A}_2\,$-moves, 
as described in {\em 4.2)}, and 
then taking the class of that link in $\,{\cal L}_M\,$. To see that  $\,P\,$ 
is well defined by this prescription, we first need to show that the image in 
$\,{\cal L}_M\,$ does not depend on the way we apply the $\,{\cal A}_2\,$-moves,
i.e., the choice of recombination ribbons, and is also invariant under the 
equivalences in $\,{\cal BL}_M\,$.

\begin{center}\ 
\epsfxsize=5.9in
\epsfbox{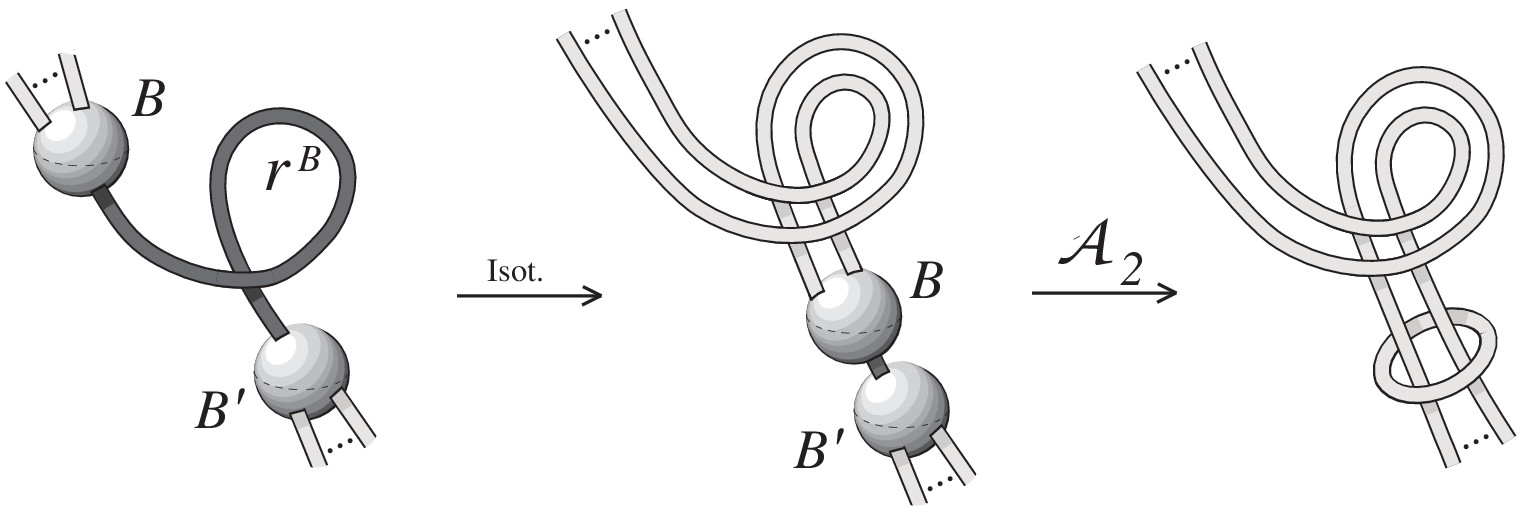}
\newline{{\sc Figure 8:} Recombination along $\,r^B\,$}
\end{center}

Consider for a pair of balls in a bridged link two choices, $r_1^B\,$ and $r_2^B\,$,
of recombination ribbons. As depicted in Figure 9 we obtain either the 
link on the left with annulus $A_2\,$ or the link on the right  with annulus $A_1\,$. 
To both diagrams we can apply an $\,{\cal O}_3\,$-move as indicated, where
the long ribbon $\,C\,$ runs along $r_1^B\,$ and $r_2^B\,$, and the short one
is the
respective opposite annulus. The resulting links differ only by 
$\,{\cal O}_2\,$-slides over $\,C\,$ as depicted in the bottom row of
Figure 9. Hence the links obtained from
either recombining along $r_1^B\,$ or along $r_2^B\,$ are in the same class in
$\,{\cal L}_M\,$.
\bigskip

\begin{center}\ 
\epsfxsize=6in
\epsfbox{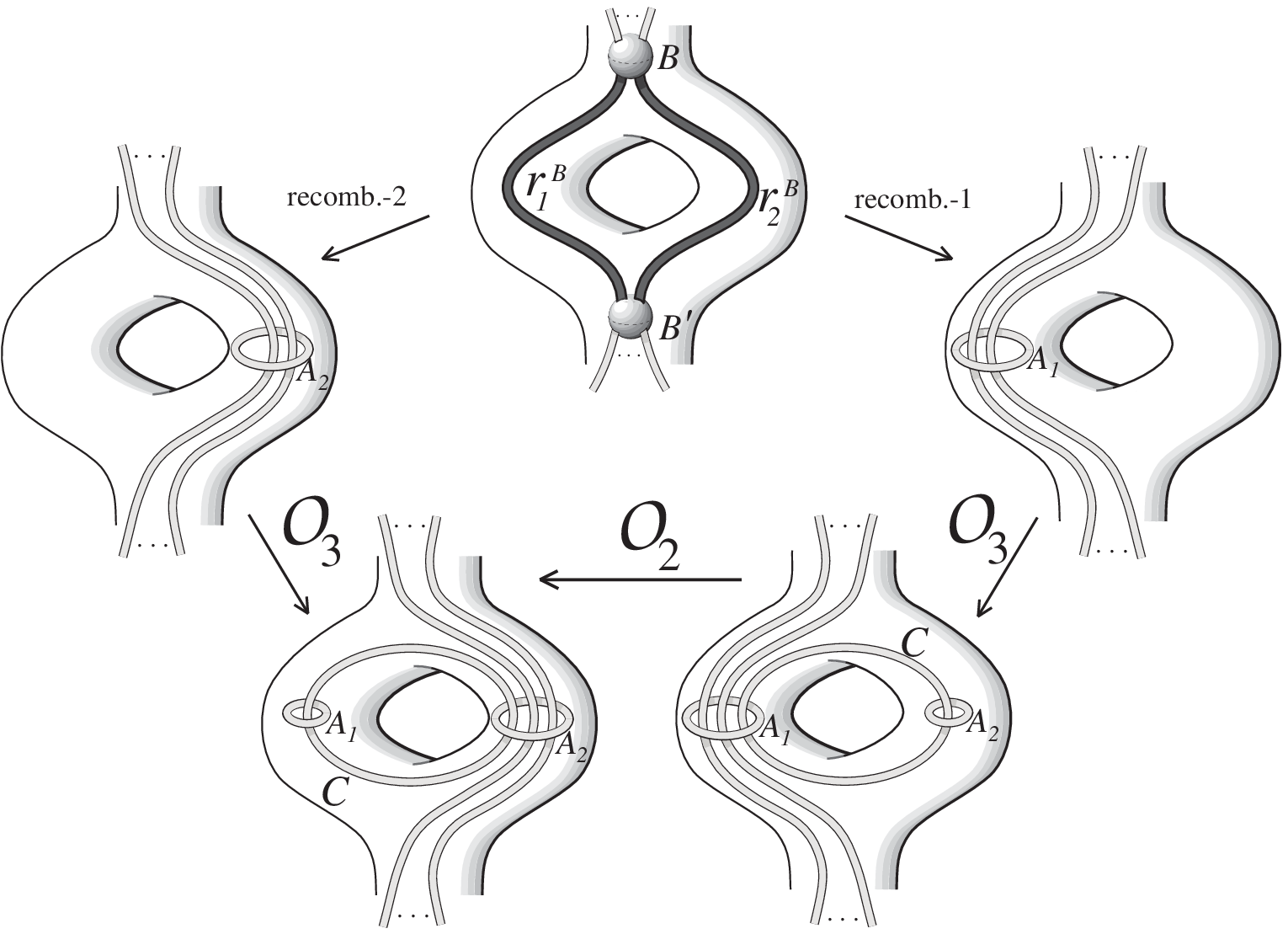}
\newline{{\sc Figure 9:} Recombinations related by  $\,{\cal O}_2\,$ {\em and} 
$\,{\cal O}_3\,$}
\end{center}

Furthermore, if two bridged links differ by an $\,{\cal A}_2\,$-move, as in 
Figure~4, and we choose the recombination ribbon to be a straight strip between
the balls, the images of the two links in Figure~4 under $\,P\,$ are 
tautologically the same. 

Finally, suppose two bridged links differ by an  $\,{\cal A}_1\,$-move, i.e.,
by an isolated configuration as in Figure 3. This configuration is mapped by
$\,P\,$ to the Hopf link $^0\bigcirc\mkern-13mu\bigcirc^0\,$ for an appropriately
chosen recombination. As described in Paragraph 1 removing or adding this link 
is a special case of the $\,{\cal O}_3\,$-move. 
\smallskip

 We can thus conclude that $\,P\,$ is well defined on $\,{\cal BL}_M\,$. It is
also obvious that $\,P\circ I\,$ is the identity in  $\,{\cal L}_M\,$ so that
$\,I\,$ must be injective. Thus $\,I\,$ is a bijection as asserted in Theorem 1.

\ \hfill{\rule{3mm}{3mm}}
\bigskip

\paragraph*{Acknowledgments:} I thank 
Robion Kirby and Justin Roberts for interest, and encouragement, and 
the referee for comments that improved  this letter in content and style.  
The author was partially supported by NSF grant DMS-9305715.

\newpage

\paragraph*{References :}\

\begin{enumerate}

\item[[FR]] Fenn, R., Rourke, C.: {On Kirby's Calculus of Links},
{\em Topology}, Vol.  {\bf 18} (1979), 1-15.

\item[[Ke]] Kerler, T.: Bridged Links and Tangle Presentations
of Cobordism Categories. \newline
 {\em Adv. Math.} to appear.(Preprint 1994).\newline
Available at {\tt http://www.math.ohio-state.edu/$\sim$kerler/papers/BL/}.

\item[[Ki]] Kirby, R.: A Calculus for framed Links in $S^3\,$,
{\em Invent. Math.} {\bf 45} (1978), 35-56.

\item[[KL]] Kerler, T., Lyubashenko, V.: Non-semisimple, Extended TQFT's.
   Preprint.

\item[[LW]] Lickorish, W.B.R.: A Representation of Orientable 3-Manifolds.
 {\em Ann. Math.} {\bf 76} (1962) 531-540.

Wallace, A.H.: Modifications and Cobounding Manifolds, {\em Cam. J. Math.}
{\bf 12 } (1960), 503-528.

\item[[R]] Roberts, J.: Kirby Calculus in Manifolds with Boundary. Preprint 1996.

\end{enumerate}

\bigskip

\bigskip

{\small {\sc Harvard University, Cambridge, MA, USA}

present address:
 
{\small {\sc The Ohio State University, Columbus,  OH, USA}

{\em E-mail-address:} kerler@math.ohio-state.edu }

\end{document}